\documentclass[3p,10pt,twocolumn]{elsarticle}

\usepackage{amssymb}

\newcounter{theorem}
\newtheorem{theorem}{Theorem}
\newcounter{lemma}
\newtheorem{lemma}{Lemma}
\newcounter{corollary}
\newtheorem{corollary}{Corollary}
\newcounter{example}
\newtheorem{example}{Example}
\newcounter{proposition}
\newtheorem{proposition}{Proposition}
\newcommand{\lbrac}{\{\!\!\{}
\newcommand{\rbrac}{\}\!\!\}}

\begin{document}

\begin{frontmatter}


    
\title{Rational complex B\'ezier curves}


\author{A. Cant\'on}
\author{L. Fern\'andez-Jambrina}
\author{M.J. V\'azquez-Gallo}
\address{Matem\'atica e Inform\'atica Aplicadas a las Ingenier\'\i as 
Civil y Naval\\
Universidad Polit\'ecnica de Madrid\\
E-28040-Madrid, Spain}

\begin{abstract}
In this paper we develop the formalism of rational complex B\'ezier
curves.  This framework is a simple extension of the CAD paradigm,
since it describes arc of curves in terms of control polygons and
weights, which are extended to complex values.  One of the major
advantages of this extension is that we may make use of two different
groups of projective transformations.  Besides the group of projective
transformations of the real plane, we have the group of complex
projective transformations.  This allows us to apply useful
transformations like the geometric inversion to curves in design.  In
addition to this, the use of the complex formulation allows to lower
the degree of the curves in some cases.  This can be checked using the
resultant of two polynomials and provides a simple formula for
determining whether a rational cubic curve is a conic or not.
Examples of application of the formalism to classical curves are
included.

\end{abstract}


\end{frontmatter}


\section{Introduction\label{intro}}


Rational spline parametrisations are the usual form for representing curves 
in CAGD \cite{farin}. In particular, using Bernstein polynomials 
\[B^n_{j}(t)=\frac{n!}{j! (n-j)!}t^{j}(1-t)^{n-j},\quad j=0,\ldots,n\] curves of degree $n$ with just one piece are parametrised 
as B\'ezier curves,
\begin{equation}\label{param}
c(t)=\frac{\displaystyle\sum_{j=0}^n\omega_{j}c_{j}B^n_{j}(t)}
{\displaystyle\sum_{j=0}^n\omega_{j}B^n_{j}(t)},\quad t\in[0,1],\end{equation}
in terms of a list of points $ \{ c_{0},\ldots,c_{n}\}$ of $\mathbb{R}^p$  named 
\emph{control polygon} and a list of coefficients 
$\{\omega_{0},\ldots,\omega_{n}\}$ dubbed \emph{weights}.

Alternatively, the rational curve in $\mathbb{R}^p$ may be viewed as 
a polynomial curve in $\mathbb{R}^{p+1}$,
\begin{equation}\label{vectorial}
\mathbf{c}(t)=\sum_{j=0}^n\mathbf{c}_{j}B^n_{j}(t),\quad t\in[0,1],
\end{equation}
using a control polygon $ \{ \mathbf{c}_{0},\ldots,\mathbf{c}_{n}\}$ formed by vectors 
$\mathbf{c}_{j}=(\omega_{j},\omega_{j}c_{j})$.

We can switch from the polynomial to the rational form dividing by the
first component of (\ref{vectorial}) and projecting into
$\mathbb{R}^p$. And conversely, the denominator of (\ref{param}) 
provides  the control vertices with the extra coordinate of the parametrisation in (\ref{vectorial}).


Rational B\'ezier parametrisations are barycentric combinations of the
control points, $\{c_{0},\ldots,c_{n}\}$, since the sum of their 
coefficients in (\ref{param}) amounts to unity.

This property allows rational parametrisations to inherit nice properties from B\'ezier 
curves, such as the convex hull property or invariance under affine 
applications \cite{farin}. 

If we deal just with planar curves, we may view points in the plane as
complex numbers and we may resort to complex geometry, adding new
features to the previous properties and results, since we may take 
advantage of both real geometry for $\mathbb{R}^{2}$ and complex 
geometry.  

Weights are positive real numbers, but in a complex framework their definition
may be extended to complex numbers.  For instance, in \cite{farouki}
Pythagorean-hodograph curves are used for computing rational offsets
to rational curves.  In \cite{leonardo-chord}, complex curves are used
for characterising chord-length parametrisations.  This formalism was
further developed in \cite{reyes-complex} using Farin or weight
points.  In \cite{conics} complex geometry is used for determining
geometric characteristics of conics.  This work was extended to
quadrics in \cite{oldquadric, quadric}.  A complex de Casteljau
algorithm is devised in \cite{schicho}.

The aim of this paper is a further development of the complex
formalism for rational planar curves.  With this purpose, it is
organised in the following fashion.  After this introduction, we
define in the next section rational complex B\'ezier curves and show
their properties.  As an example of this construction, in
Section~\ref{arcs} we produce arcs of circles.  In order to relate
real and complex parametrisations, in Section~\ref{irred} we introduce
the concept of reducible parametrisations and provide
characterisations for it and we show in Section~\ref{factors} how to
simplify them in B\'ezier form.  Having shown how to derive complex
parametrisations from real ones, we go the way back in
Section~\ref{real}.  Section~\ref{curves} is devoted to examples of
curves which can be parametrised as complex rational curves using
inversion.  Finally in the last section we derive some conclusions.

\section{Rational complex B\'ezier curves\label{define}}

We may identify $\mathbb{R}^2$ with the complex plane $\mathbb{C}$ through
the usual bijection
\[\begin{array}{ccc}\mathbb{R}^2&\rightarrow&\mathbb{C}\\ 
(x,y)&\mapsto&x+iy,\end{array}\] 
so that control polygons of planar curves are readily written in complex form 
$\{z_{0},\ldots,z_{n}\}$, with $z_{k}={c_{k}^x}+i{c_{k}^y}$, 
$k=0,\ldots,n$, for $c_{k}=({c_{k}}^{x},{c_{k}}^{y})$.


A rational complex B\'ezier curve of degree $n$ is then a planar curve
which can be parametrised as
\begin{equation}\label{cparam}
c(t)=\frac{\displaystyle \sum_{j=0}^n\omega_{j}z_{j}B^n_{j}(t)}
{\displaystyle \sum_{j=0}^n\omega_{j}B^n_{j}(t)} ,\quad 
t\in[0,1],
\end{equation}
in terms of the Bernstein polynomials of degree $n$ 
and two 
lists of complex numbers: the complex
\emph{control polygon} $\{z_{0},\ldots,z_{n}\}$ and the complex
\emph{weights} $\{\omega_{0},\ldots,\omega_{n}\}$ of the curve, where 
at least $\omega_{0}$ and $\omega_{n}$ are non-zero.

Most properties of real B\'ezier curves \cite{farin} are inherited by complex 
curves:

\begin{itemize}
    \item  $z_{0},z_{n}$ are the endpoints of the curve: 
    \[c(0)=z_{0},\qquad c(1)=z_{n}.\]

    \item  The tangent vectors at the endpoints of the curve,
\[\hspace{-0.7cm}c'(0)=\frac{n\omega_{1}}{\omega_{0}}(z_{1}-z_{0}),\ 
c'(1)=\frac{n\omega_{n-1}}{\omega_{n}}(z_{n}-z_{n-1}),\]
consist of a rotation and a stretching of the first and last edges of
the control polygon.  

For instance, if $\omega_{1}/\omega_{0}=r
e^{i\phi}$, the tangent vector is a modification of the vector
$\overrightarrow{z_0z_1}$ by a factor $nr$ and a rotation of an angle
$\phi$. And analogously at the other endpoint of the curve, 
where the stretching and rotation are given by 
$n\omega_{n-1}/\omega_{n}$.

\item Projective invariance: The affine invariance property of 
B\'ezier curves can be
extended to projective invariance.  Any regular linear transformation
of $\mathbb{R}^{p+1}$ defines a
projective transformation for $\mathbb{P}^p$ up to a non-zero constant.  That
is, two linear transformations of $\mathbb{R}^{p+1},$ $\lambda
\mathbf{f}$ and $\mathbf{f}$, provide the same projective
transformation $f$ of $\mathbb{P}^{p}$ for 
$\lambda\in\mathbb{R}\backslash\{0\}$.

If we apply a projective
transformation $f$ to the curve, with $\mathbf{f}$ as associated
linear map,
\[ \mathbf{f}(\mathbf{c}(t))=\sum_{j=0}^n\mathbf{f}(\mathbf{c}_{j})B^n_{j}(t),\]
the transformed curve is another curve of degree $n$ which has 
$\{\mathbf{f}(\mathbf{c}_{0}),\ldots,\mathbf{f}(\mathbf{c}_{n})\}$ as vector control 
polygon  \cite{farinNURBS}.

That is, we have the invariance under real projective 
transformations,
\begin{equation}\label{realproy}\hspace{-0.7cm}
f(x,y)=\left(\frac{a_{1}+b_{1}x+c_{1}y}{a_{0}+b_{0}x+c_{0}y},
\frac{a_{2}+b_{2}x+c_{2}y}{a_{0}+b_{0}x+c_{0}y}\right),\end{equation}
associated to regular linear transformations of $\mathbb{R}^3$, 
\[\mathbf{f}\left(\begin{array}{c}x_{0}\\x_{1}\\x_{2}
\end{array}\right)=
\left(\begin{array}{ccc}a_{0} & b_{0} & c_{0} \\
a_{1} & b_{1} & c_{1} \\ a_{2} & b_{2} & c_{2} 
\end{array}\right)\left(\begin{array}{c}x_{0}\\x_{1}\\x_{2}\end{array}\right).\quad
\]

But besides this general projective invariance property, based on constructing
projective transformations through linear transformations of
$\mathbb{R}^3$, we have another possibility in the plane
$\mathbb{R}^{2}$.

Rational complex parametrisations 
can be also viewed as polynomial curves in $\mathbb{C}^2$,
\[\mathbf{c}(t)=\sum_{j=0}^n\mathbf{z}_{j}B^n_{j}(t),\] where
$\mathbf{z}_{j}=(\omega_{j},\omega_{j}z_{j})$, for  $j=0,\ldots,n$.

Regular linear transformations of $\mathbb{C}^2$,
\[\hspace{-0.7cm}\mathbf{f}(z_{0},z_{1})=(az_{0}+bz_{1},cz_{0}+dz_{1}), \] 
where $a,b,c,d\in\mathbb{C}$,  $ad-bc\neq0$,
provide  complex projective transformations of the projective complex 
plane, named M\"obius transformations,
\begin{equation}\label{mobius}
f(z)=\frac{c+dz}{a+bz}.
\end{equation}


We have hence for these curves another property of invariance 
under projective transformations,
\begin{equation}
\textbf{f}(\mathbf{c}(t))=\sum_{j=0}^n\textbf{f}(\mathbf{z}_{j})B^n_{j}(t).\end{equation}

The group of real projective transformations depends on a larger number of parameters than the 
group of complex projective transformations, but it does not comprise 
it.

For instance, the inversion transformation $f(z)=1/z$,
\[f(z)=\frac{1}{x+iy}=\frac{x}{x^2+y^2}-\frac{iy}{x^2+y^2},\]
corresponding to a M\"obius transformation with $a=0=d$, $b=c$,
does not provide a real projective transformation of $\mathbb{R}^2$ 
of the form (\ref{realproy}).

Hence both groups are different and none of them is a subgroup of the 
other.

\item If $f(z)=1/z$ then $f(c(t))$ is a rational complex curve with
control polygon $\{\frac{1}{z_0},\dots,\frac{1}{z_n}\}$ and complex
weights $\{\omega_0z_0,\dots,\omega_nz_n\}$.

    \item  The parametrisation does not change if we multiply all 
    weights by a non-zero common factor 
	$\lambda\in\mathbb{C}\backslash\{0\}$, 
    $\{\lambda\omega_{0},\ldots,\lambda\omega_{n}\}$.

    \item  Another rational parametrisation of degree $n$ for the 
    same curve is obtained by transforming the parameter with a 
   real M\"obius transformation of the interval $[0,1]$ onto itself,
\begin{equation}t(u)=\frac{u}{(1-\rho)u+\rho},\quad 
	u\in[0,1],\label{mobiuspar}\end{equation}
for a non-zero real parameter $\rho$.

The control polygon remains the same, but the new weights are 
$\hat \omega_{j}=\rho^{n-j} \omega_{j}$, $j=0,\ldots,n$.

\item Degree elevation: If we multiply both the numerator and 
    the denominator of $c(t)$ by a polynomial $p(t)$ of degree one
\[   p(t)=\alpha (1-t)+ \beta t,\quad \alpha,\beta\in\mathbb{C},\]
the parametrisation is formally of degree $n+1$ \cite{denker,hoschek}
with new control polygon $\{z^1_{0},\ldots,z_{n+1}^1\}$ and weights 
$\{\omega^1_{0},\ldots,\omega_{n+1}^1\}$, 
\begin{eqnarray}\label{elevrac}\hspace{-0.7cm}
    \omega^1_{j}\!\!\!\!&=&\!\!\!\!\alpha \frac{n+1-j}{n+1}\omega_{j}+
    \beta \frac{j}{n+1}\omega_{j-1},\nonumber\\\hspace{-0.7cm}
    \omega^1_{j}z^1_{j}\!\!\!\!&=&\!\!\!\!\alpha \frac{n+1-j}{n+1}\omega_{j}z_{j}+
    \beta \frac{j}{n+1}\omega_{j-1}z_{j-1}.
\end{eqnarray}

The possibility of having complex $\alpha$ and $\beta$ increases the 
ways of performing degree elevation compared to the real case. 

This poses the question of whether a complex parametrisation can be 
further simplified or not. We come back to  this issue in 
Sections~\ref{irred} and \ref{factors}. In fact, converting complex 
parametrisations into real ones is a case of degree-elevation, as we 
shall see in more detail in Section~\ref{real}.

\item Convex envelope property: if weights are not real and positive, 
the curve does not lie in the convex envelope of its control polygon. 

We lose this feature, but complex weights add flexibility 
for designing curves, as we see in the following section.

\end{itemize}



\section{Arcs of circles\label{arcs}}

The image of a line segment by a M\"obius transformation as in 
(\ref{mobius}) is
either a line segment or an arc of a circle when the denominator does
not vanish at any point of the segment.  

Considering M\"obius transformations that fix $z_0$ and $z_1$, 
\[f(z)=\frac{\sigma z_{1}(z-z_{0})-z_{0}(z-z_{1})}{\sigma (z-z_{0})-(z-z_{1})},\]
for $\sigma\in\mathbb{C}\backslash\{0\}$, the image of the line segment 
linking those points, $c(t)=z_0(1-t)+z_1t$, $t\in[0,1]$, is
\[f(c(t))=\frac{z_{0}(1-t)+\sigma z_{1}t}{(1-t)+\sigma t},
\]
and we may write it in B\'ezier form taking $\sigma =\omega_{1}/\omega_{0}$, 
where $\omega_0,\omega_1$ are non-zero complex numbers.  

The curve $f(c(t))$ is an arc of a circle if $\sigma $ is not a real number,
so from now on this will be assumed.

\begin{figure}
\centering
\includegraphics[height=6cm]{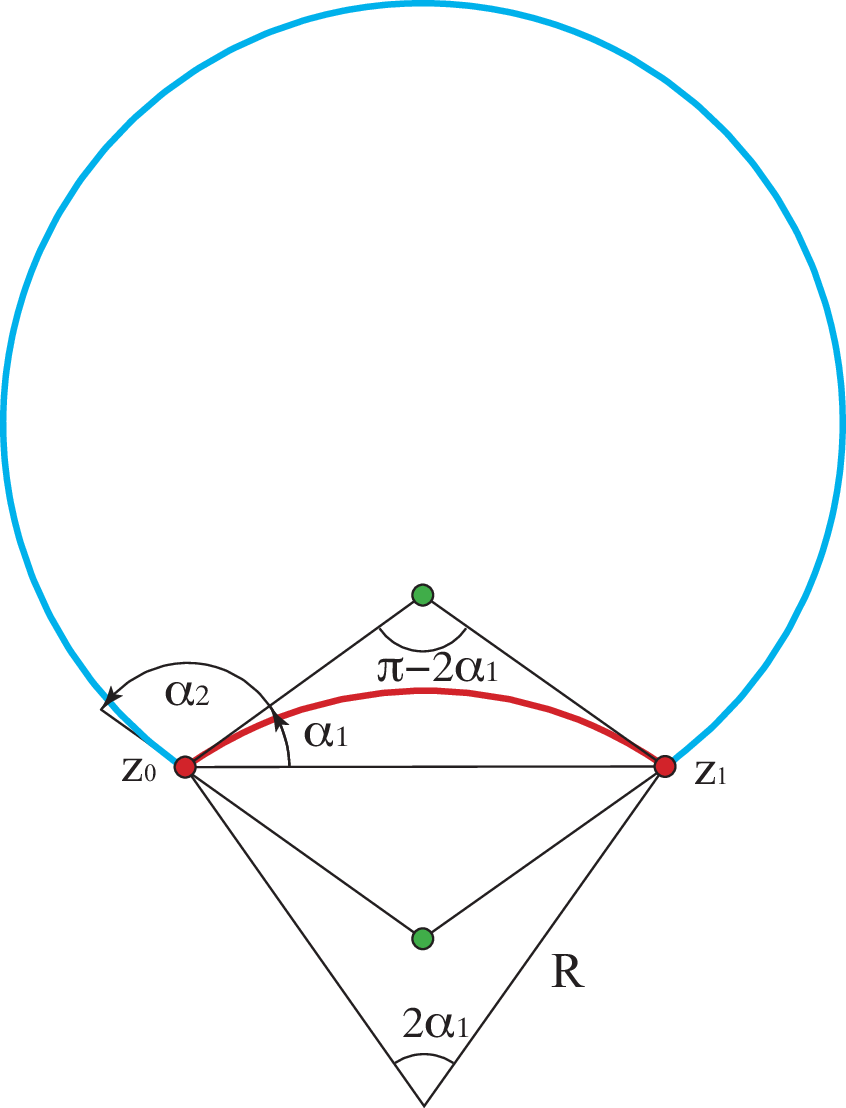}
\caption{Two arcs of circle from $z_{0}$ to $z_{1}$}
\label{circle}
\end{figure}

As mentioned in Section~\ref{define}, the velocities at both ends of the arc of the circle are
\[
c'(0)=\sigma(z_1-z_0),\quad c'(1)=\sigma^{-1}(z_1-z_0).
\]
If $\sigma=re^{i\alpha}$ with $0<|\alpha|<\pi$ and $r>0$, 
then $c'(0)$ is a rotation of angle $\alpha$ and a stretching with 
factor $r$ of the vector $z_1-z_0$. Moreover, the parameter $r$ can be obtained by choosing a
particular point on the arc of the circle to be the image by the M\"obius 
transformation of the midpoint of the segment (a M\"obius transformation is completely 
determined by fixing the image of three different points). When $r=1$, the midpoint of 
the segment is mapped to the midpoint of the arc of the circle. Thus 
performing a change of parameter as in (\ref{mobiuspar}), it can be obtained that 
$\sigma=e^{i\alpha}$ and the arcs of the circles 
in Figure~\ref{circle} can be described by the rational complex curve
\[
c(t)=\frac{z_0(1-t)+e^{i\alpha}z_1t}{(1-t)+e^{i\alpha}t}.
\]

The radius of the circle is:
\[
R=\frac{|z_1-z_0|}{2|\sin \alpha|},
\]
and its center,
\[
z_0+\frac{ie^{i\alpha}}{2\sin\alpha}(z_0-z_1).
\]

In Figure~\ref{circle} we may see two arcs, for
$\alpha_{1}=\pi/6$, $\alpha_{2}= 5\pi/6$ and the same complex control
polygon, but different weights.  The green control points belong to 
control polygons for arcs as curves of degree two, whereas as complex 
curves we just need $z_{0}$ and $z_{1}$.


Therefore, any circular arc can be represented by a complex rational
B\'ezier curve of degree one, whereas its degree is two as a real
rational B\'ezier curve.  Hence, when written in complex notation the
numerator and denominator of the rational expression have a common
factor (the rational function is reducible).  

In the next sections we
are going to show how to identify common factors in complex rational
curves written in B\'ezier form.

%
%

\section{Degree elevation and irreducibility\label{irred}}


There is a simple way of checking whether a rational complex
parametrisation of a planar curve, $c(t)=p(t)/q(t)$, cannot be
further simplified, that is, if $c(t)$ is irreducible.


A complex rational parametrisation is irreducible if and only if
the polynomials $p(t)$ and $q(t)$ are relatively prime.  This means it cannot be
further simplified by cancellation of common complex polynomial factors and 
the degree of the parametrisation is the lowest possible.

This is specially useful, since the simplest way of casting a 
rational parametrisation in complex form
\[c(t)=\left(\frac{x(t)}{w(t)},\frac{y(t)}{w(t)}\right)\rightarrow
\frac{x(t)+i y(t)}{w(t)},\]
often produces polynomials $z(t)=x(t)+i y(t)$, $w(t)$ which may be simplified by
cancellation of common polynomial factors.  As in the case of the
parametrisation of an arc of a circle, which is of degree two as a
real parametrisation, but of degree one as a complex one.


In order to know if a parametrisation is irreducible, we may use the 
Sylvester resultant $R(p,q)$ of two polynomials $p$, $q$ \cite{walker}.

If we denote by $T(c_k,\dots,c_0;l)$ the matrix with $k+l$ columns 
and $l$ rows given by 
\[\hspace{-1cm}T(c_{k},\ldots,c_{0};l):=\left(\begin{array}{ccccc}
c_{k} &\cdots &c_{0} & 0 &\cdots \\
0 &\ddots &\ddots & \ddots &\ddots \\
0 & 0 & c_{k} &\cdots &c_{0}
\end{array}\right)
\begin{array} {c}
\uparrow \\
l\ \textrm{rows}\\
\downarrow 
\end{array},\]
the Sylvester resultant $R(p,q)$ of two polynomials $p$, $q$,
\[p(t)=\sum_{j=0}^{m}a_{j}t^{j},\qquad q(t)=\sum_{j=0}^{n}b_{j}t^{j},\]
respectively of degree $m>0$ and $n>0$, is defined as
\begin{equation}\label{resultant0}
R(p,q):=\left|\begin{array}{c}T(a_{m},\ldots,a_{0};n)\\
T(b_{n},\ldots,b_{0};m)\end{array}\right|
,\end{equation}
which is a determinant with $m+n$ rows and columns.

%

The resultant of two polynomials is related to the complex roots of both 
polynomials by a classical result \cite{walker}:
\begin{theorem}Let $p(t)=\displaystyle\sum_{j=0}^{m}a_{j}t^{j}$,
$q(t)=\displaystyle\sum_{j=0}^{n}b_{j}t^{j}$ be two polynomials with respective 
complex roots $\{x_{1},\ldots,x_{m}\}$, $\{y_{1},\ldots,y_{n}\}$, then
\[
R(p,q)=(a_{m})^{n}(b_{n})^{m}\prod_{1\le j\le m}\prod_{1\le k\le 
n}(x_{j}-y_{k}).\]
\end{theorem}

As a consequence, the resultant of two polynomials vanishes if they 
share at least one root $x_{J}=y_{K}$. In this case, they share a 
common factor $t-x_{J}$ and their quotient is reducible:

\begin{corollary}The quotient of two polynomials $p(t)$ and 
$q(t)$, $c(t)=p(t)/q(t)$ is irreducible if and only if $R(p,q)\neq0$.
\end{corollary}

The resultant of two polynomials is written in (\ref{resultant0}) in
terms of their coefficients in the monomial basis
$\{1,t,t^{2},\ldots\}$.  Since we are interested in applying it to
CAGD, it is useful to write it in the Bernstein basis, and
this has already been done in \cite{winkler}:  
\begin{theorem}The resultant of
two polynomials $p(t)=\displaystyle\sum_{j=0}^{m}a_{j}B^{m}_{j}(t)$,
$q(t)=\displaystyle\sum_{j=0}^{n}b_{j}B^n_{j}(t)$ in Bernstein form
may be written as
\begin{equation}\hspace{-0.5cm}
R(p,q):=D\left|\begin{array}{c}T
\left(a_{0}{m\choose{0}} ,\ldots,a_{m}{m\choose{m}} ;n\right)\\\\
T\left(b_{0}{n\choose{0}} ,\ldots,b_{n}{n\choose{n}} ;m\right)\end{array}\right|
,\label{resultant}
\end{equation}
\[D^{-1}=\prod_{j=0}^{m+n-1}{{m+n-1}\choose{j}},\]
\label{thresultant}\end{theorem}
and this provides this straightforward result:
\begin{corollary}The quotient $c(t)=p(t)/q(t)$ of two polynomials
$p(t)=\displaystyle\sum_{j=0}^{m}a_{j}B^{m}_{j}(t)$,
$q(t)=\displaystyle\sum_{j=0}^{n}b_{j}B^n_{j}(t)$ in Bernstein form
is irreducible if and only if
\[
\left|\begin{array}{c}T
\left(a_{0}{m\choose{0}} ,\ldots,a_{m}{m\choose{m}} ;n\right)\\\\
T\left(b_{0}{n\choose{0}} ,\ldots,b_{n}{n\choose{n}} ;m\right)\end{array}\right|
\neq0.\]
\end{corollary}

As we see the condition of non-vanishing resultant is pretty similar
in both bases, except for the use of coefficients $\tilde
a_{j}:=a_{j}{m\choose{j}}$, $\tilde b_{k}:=b_{k}{n\choose{k}}$ instead
of $a_{j}$, $b_{k}$.  

We shall call \emph{reduced} the coefficients with the tilde.  In
order to tell reduced lists from standard one, we shall use double
braces $\lbrac, \rbrac$ for the reduced ones if necessary to avoid
confusions.  Of course, for polynomials of degree one, reduced
coefficients are the same as the standard ones.

Using this reduced notation, the resultant is written as
\[
R(p,q)=D\left|\begin{array}{c}T(\tilde a_{0},\ldots,\tilde a_{m};n)\\
T(\tilde b_{0},\ldots,\tilde b_{n};m)\end{array}\right|.\]

In our case, we have a rational parametrisation (\ref{cparam}) of a
curve as a quotient of complex polynomials, being 
$\{\omega_{0},\ldots,\omega_{n}\}$ the coefficients of its 
denominator and $\{\omega_{0}z_{0},\ldots,\omega_{n}z_{n}\}$ the coefficients 
of its numerator. 

Applying the previous corollary to these polynomials, we may write a 
characterisation of irreducibility of complex
B\'ezier parametrisations in terms of their weights and control
polygons:
\begin{corollary}
A rational B\'ezier parametrisation with complex control polygon 
$\{z_0,\ldots,z_{n}\}$ and complex weights
$\{\omega_{0},\ldots,\omega_{n}\}$ is
irreducible if and only if
\[
\left|\begin{array}{c}
T\left(\tilde\omega_{0}z_{0},\ldots,\tilde\omega_{n}z_{n};n\right)\\
T\left(\tilde\omega_{0},\ldots,\tilde\omega_{n};n\right)

\end{array}\right|
%
\neq0.
\]
\end{corollary}

In this corollary, for the sake of simplicity we have included the
Bernstein polynomials in what we have called the reduced version of the
weights $\{\tilde\omega_0,\ldots,\tilde\omega_{n}\}$.

\begin{example}\label{excircle} The arc of a circle described by the real 
control polygon $\{(1,0),(1,1),(0,1)\}$ and weights 
$\{1/2,1/2,1\}$, which produces the usual parametrisation
\[c(t)=\left(\frac{1-t^{2}}{1+t^{2}},\frac{2t}{1+t^{2}}\right),\quad 
t\in[0,1],\]
can be seen to be reducible.

The complex control polygon $\{1,1+i,i\}$ and its weighted 
complex control polygon $\{1/2,(1+i)/2,i\}$,
\[\hspace{-1cm}
\left|\begin{array}{cccc}\omega_{0}z_{0}&2\omega_{1}z_{1}&\omega_{2}z_{2}&0\\
0&\omega_{0}z_{0}&2\omega_{1}z_{1}&\omega_{2}z_{2}\\
\omega_{0}&2\omega_{1}&\omega_{2}&0\\
0&\omega_{0}&2\omega_{1}&\omega_{2}\end{array}\right|=
\left|\begin{array}{cccc}\frac{1}{2}&1+i&i&0\\
0&\frac{1}{2}&1+i&i\\
\frac{1}{2}&1&1&0\\
0&\frac{1}{2}&1&1\end{array}\right|=0.\]
are seen to produce a vanishing resultant.\end{example}


The case of conic sections which are written as rational cubic 
B\'ezier curves has attracted some attention \cite{wang-conic, conicubicas}. In this case rational 
parametrisations of degree two are expressed as formal cubics through 
degree elevation \cite{denker}, which involves multiplying numerator and 
denominator of the parametrisation by a linear polynomial. In such 
case, the parametrisation is reducible and we have:
\begin{theorem}
A rational B\'ezier curve described by its complex control polygon 
$\{z_{0},z_{1},z_{2},z_{3}\}$ and weights 
$\{\omega_{0},\omega_{1},\omega_{2},\omega_{3}\}$ is a conic arc if and only if
\[\hspace{-0.8cm}\left|\begin{array}{cccccc}
\omega_{0}z_{0}&3\omega_{1}z_{1}&3\omega_{2}z_{2}&\omega_{3}z_{3}&0&0\\
0&\omega_{0}z_{0}&3\omega_{1}z_{1}&3\omega_{2}z_{2}&\omega_{3}z_{3}&0\\
0&0&\omega_{0}z_{0}&3\omega_{1}z_{1}&3\omega_{2}z_{2}&\omega_{3}z_{3}\\
\omega_{0}&3\omega_{1}&3\omega_{2}&\omega_{3}&0&0\\
0&\omega_{0}&3\omega_{1}&3\omega_{2}&\omega_{3}&0\\
0&0&\omega_{0}&3\omega_{1}&3\omega_{2}&\omega_{3}
\end{array}\right|=0.\]
\end{theorem}

This theorem is just Corollary 3 in the case of cubic 
parametrisations and it allows characterisation of cubic parametrisations of conic arcs in
terms of a complex condition, which is equivalent to the two
conditions in \cite{wang-conic}.

\begin{example}The cubic B\'ezier curve with control polygon $\{(1,
0), (1, 4/5), (1/2, 1), (0, 1)\}$ and weights $\{2, 5/3, 4/3,1\}$ is in 
fact an arc of a conic.\label{fakecubic}

The complex control polygon for the curve is $\{1, 1+4i/5,1/2+i,i\}$.

We may compute the Sylvester resultant of the parametrisation,
\[\left|\begin{array}{cccccc}
2&5+4i&2+4i&i&0&0\\
0&2&5+4i&2+4i&i&0\\
0&0&2&5+4i&2+4i&i\\
2&5&4&1&0&0\\
0&2&5&4&1&0\\
0&0&2&5&4&1
\end{array}\right|\] and check that it vanishes, so that the 
parametrisation corresponds in fact to a conic arc.\end{example}

Finally, multiplication of polynomials in Bernstein form is simpler
and very similar to multiplication in monomial basis, using the
reduced form, just taking into account the powers of $(1-t)$ and $t$. 
We shall make use of it in the following section.

If we multiply two polynomials with reduced coefficients
$\{\tilde a_{0},\ldots, \tilde a_{m}\}$ and $\{\tilde
b_{0},\ldots, \tilde b_{n}\}$, we get as a resulting list of reduced
coefficients for the product,\begin{eqnarray}\{\tilde 
c_{0},\ldots,\tilde c_{m+n}\}
\{\tilde a_{0},\ldots, \tilde a_{m}\}* \{\tilde
b_{0},\ldots, \tilde b_{n}\},\nonumber\\
\tilde c_{J}=\sum_{j+k=J}\tilde a_{j}\tilde 
b_{k},\ J=0,\ldots,m+n,\label{multiply}\end{eqnarray}
which is the convolution of both lists.


\section{Factorisation\label{factors}}

Once we know a complex parametrisation is reducible, it would be nice
to have a procedure for factoring out the common terms in numerator
and denominator, both of them written in the Bernstein basis.  The
common factor is the greatest common divisor of both polynomials.  We
would also need a procedure for dividing polynomials in the Bernstein
basis.

For singling out the greatest common divisor of two polynomials, we 
may use Euclid's algorithm for polynomials in Bernstein's basis.

In \cite{buse} there are versions of the algorithms for
dividing polynomials and computing their greatest common divisor in
their Bernstein form:
\begin{proposition}\label{propdiv}
If $p(t)=\displaystyle\sum_{j=0}^{m}a_{j}B^{m}_{j}(t)$ and
$q(t)=\displaystyle\sum_{j=0}^{n}b_{j}B^{n}_{j}(t)$ are two
polynomials with $m\ge n\ge0$, such that $a_m=p(1)\neq 0$,
$b_n=q(1)\neq 0$, then there exist two uniquely determined polynomials
$f(t)$ (quotient) and $r(t)$ (remainder), respectively of degrees
$m-n$ and $n-1$, such that
\[
p(t)=q(t)f(t)+(1-t)^{m-n+1}r(t).
\]
\end{proposition}

Notice that this division algorithm provides an actual \emph{remainder} of degree 
$m$, since, though $r(t)$ has formally degree $n-1$, it is multiplied by 
a factor of degree $m-n+1$. In this sense, the algorithm just 
eliminates the contribution of polynomials $t^{m}$, $t^{m-1}(1-t)$,\ldots, 
$t^{n}(1-t)^{m-n}$. 

The division algorithm works even if $a_m= 0$, but we have to factor
powers of $t$ beforehand.

Another division algorithm can be produced for eliminating the 
contribution of high powers of $(1-t)$ instead, requiring 
$a_{0}=p(0)\neq0$, $b_{0}=q(0)\neq0$.

These requirements on the polynomials are not a big issue, since common $(1-t)$ and $t$ factors are
easily introduced and factored. 
\begin{lemma}If a polynomial $p(t)=\sum_{j=0}^{n}a_{j}B^{n}_{j}(t)$ of degree $n$ has
reduced coefficients $\{\tilde a_{0},\ldots, \tilde a_{n}\}$, then the
polynomial $q(t)=(1-t)^{J}p(t)t^{K}$, has reduced coefficients
$\{\underbrace{0,\ldots, 0}_{K\ \mathrm{times}},\tilde 
a_{0},\ldots,
\tilde a_{n}, \underbrace{0,\ldots, 0}_{J\
\mathrm{times}}\}$.


If a polynomial $p(t)$ of degree $n$ has
reduced coefficients $\{0,\ldots,0, \tilde a_{J},\ldots, \tilde
a_{K},0,\ldots,0\}$, $0\le J<K\le n$, then it can be expressed as 
$p(t)=(1-t)^{n-K}q(t)
t^{J}$, where the polynomial $q(t)$ has reduced coefficients $\{\tilde
a_{J},\ldots ,\tilde a_{K}\}$.\end{lemma}

This lemma is readily proven, taking into account the missing or
added powers of $t$ and $(1-t)$.  The use of lists of reduced coefficients
is useful for avoiding 
cumbersome binomial coefficients on changing the degree of
polynomials.
%
%

\begin{example}The polynomial $p(t)=2(1-t)^{2}t^{2}+(1-t)t^{3}$ with
reduced coefficients $\lbrac0,0,2,1, 0\rbrac$.

It can be written as $p(t)=(1-t)t^{2}q(t)$, $q(t)=2(1-t)+t$, where 
the coefficients of $q(t)$ are $\lbrac2,1\rbrac$.

\begin{example}Divide $p(t)=4t^4 + t^3 - t^2 - t + 1$ by 
$q(t)=-t^3 + 6t^2 - 5t + 2$ in Bernstein form.\end{example}

The coefficients in the Bernstein basis are seen to be $\{1, 3/4, 1/3, 0, 4\}$ ($\lbrac1,3,2,0,4\rbrac$ in reduced 
form) for $p(t)$ and $\{2, 1/3, 2/3, 2\}$ ($\lbrac2,1,2,2\rbrac$ in reduced form) for $q(t)$.

Then $p(t)=q(t)f(t) +(1-t)^{2}r(t)$, where the quotient $f(t)$ is of 
degree one and $r(t)$ is of degree two. Their reduced coefficients 
would be $\{\tilde f_{0},\tilde f_{1}\}$ and 
$\{\tilde r_{0},\tilde r_{1}, \tilde r_{2}\}$.

The reduced coefficients of $q(t)f(t)$ are
\[\{2\tilde f_{0},\tilde f_{0}+2\tilde f_{1},2\tilde f_{0}+
\tilde f_{1},2\tilde f_{0}+2\tilde f_{1},2\tilde f_{1}\},\]
and we get for $p(t)-f(t)q(t)$,
\[\hspace{-1cm}\{1-2\tilde f_{0},3-\tilde f_{0}-2\tilde f_{1},2-2\tilde f_{0}-
\tilde f_{1},-2\tilde f_{0}-2\tilde f_{1},4-2\tilde f_{1}\},\]
from which we learn $\tilde f_{1}=2$, $\tilde f_{0}=-2$.

For the remainder we get from the reduced coefficients of $p(t)-f(t)q(t)$,
$\lbrac5,1,4,0,0\rbrac$, and after eliminating zeros, $\lbrac5,1,4\rbrac$.

Hence, the coefficients in the Bernstein basis of $f(t)$ and $r(t)$ are respectively 
$\{-2,2\}$, $\{5,1/2,4\}$, corresponding to polynomials,
\[f(t)=4t-2,\quad r(t)=8t^2 - 9t + 5.\]\end{example}


%
%

%
Since we have a division algorithm, we have an Euclid algorithm for
calculating the greatest common divisor (gcd) of two polynomials.
Even though the division algorithm is not the standard one and does
not provide the same quotient and remainder, Euclid's algorithm 
gives out the same gcd up to a constant factor \cite{buse}.

\begin{theorem}\textbf{Euclid's algorithm:} Let $p(t)$, $q(t)$ be two 
polynomials, where the degree of $p(t)$ is equal or greater than the degree 
of $q(t)$. The following algorithm computes 
$\mathrm{gcd}\,(p(t),q(t))$:\end{theorem}
\begin{itemize}
	\item  Step 0: Divide $p(t)$ by $q(t)$ and get
	remainder $r_{1}(t)$. If $r_{1}(t)\equiv0$, then $q(t)=\mathrm{gcd}\,(p(t),q(t))$.
	
	\item  Step 1:  Divide $q(t)$ by $r_{1}(t)$ and get	remainder $r_{2}(t)$. If $r_{2}(t)\equiv0$, then 
	$r_{1}(t)=\mathrm{gcd}\,(p(t),q(t))$.

	\item  Step 2:  Divide $r_{1}(t)$ by $r_{2}(t)$ and get remainder $r_{3}(t)$. If $r_{3}(t)\equiv0$, then 
	$r_{2}(t)=\mathrm{gcd}\,(p(t),q(t))$.

	\item  \ldots

	\item  Step $N$: Divide $r_{N-1}(t)$ by $r_{N}(t)$ and get 
	remainder $r_{N+1}(t)$. If $r_{N+1}(t)\equiv0$, 
	then $r_{N}(t)=\mathrm{gcd}\,(p(t),q(t))$.

\end{itemize}

We may review now some of the previous examples:

\begin{example} Rational curve with control polygon 
$\{(1,0),(1,1),(0,1)\}$ and weights $\{1,1,2\}$.
\begin{figure}
\centering
\includegraphics[height=4cm]{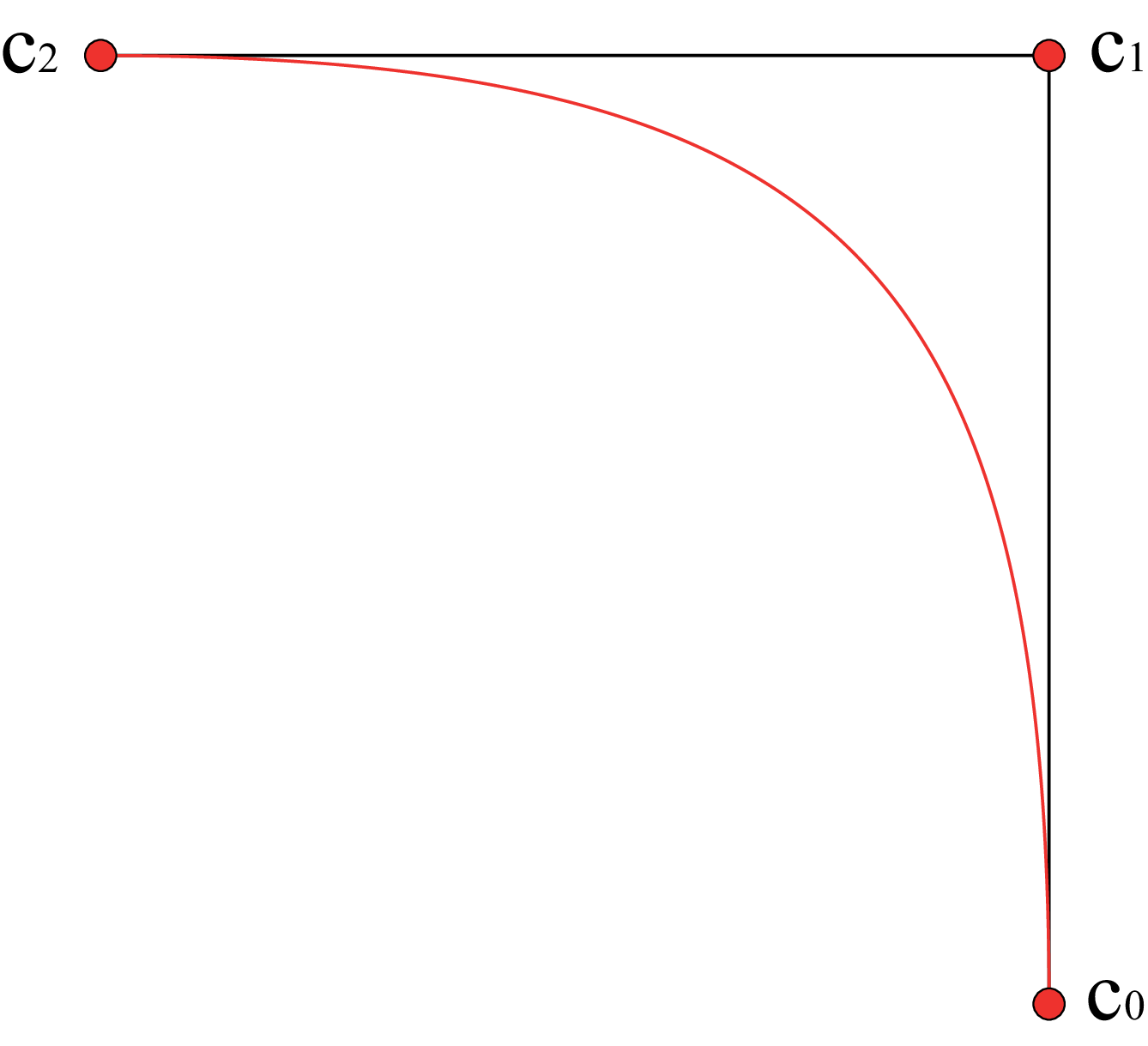}
\caption{Arc of circle with control polygon $\{1,1+i,i\}$ and 
weights $\{1,1,2\}$}
\label{circle1}
\end{figure}

We have seen this arc in Example~\ref{excircle}. It is a quarter of a circle of unit radius, extending from $(1,0)$ to 
$(0,1)$ (see Figure~\ref{circle1}):

As a complex curve, the control polygon is $\{1,1+i,i\}$, with list 
of weights $\{1,1,2\}$. This gives $\{1,1+i,2i\}$ as coefficients 
in the Bernstein basis for the numerator of the parametrisation.

We divide the polynomials $p(t)$, $q(t)$ with respective reduced 
coefficients 
$\lbrac1,2+2i,2i\rbrac$ and $\lbrac1,2,2\rbrac$ to obtain a quotient with
coefficients $\{i\}$ and $\{1-i,2\}$ for the remainder $r_{1}(t)$.


Dividing $q(t)$ by $r_{1}(t)$ we get a null remainder. 
Hence, 
$\mathrm{gcd}\,(p(t),r(t))=r_{1}(t)$.

Finally, we factor $r_{1}(t)$ from both $p(t)$ and $q(t)$: 
Dividing $q(t)$  by $r_{1}(t)$, 
we get coefficients 
$\{(1+i)/2,1\}$ for the denominator, that is, the list of weights for 
the curve.

Dividing $p(t)$  by $r_{1}(t)$, we get coefficients
$\{(1+i)/2,i\}$ for the numerator.

After factoring the weights, we 
get $\{1,i\}$ as complex control polygon for the curve.

We have already divided $q(t)$ by $r_{1}(t)$ and hence we have the
complex representation of the curve with control polygon $\{1,i\}$ and
complex weights $\{(1+i)/2,1\}$. \end{example}

We may review now the case of the degree-elevated arc of conic in 
Example~\ref{fakecubic}.

\begin{example}Cubic B\'ezier curve with control polygon $\{(1,
0), (1, 4/5), (1/2, 1), (0, 1)\}$ and weights $\{2, 5/3, 4/3,1\}$.
\begin{figure}
\centering
\includegraphics[height=4cm]{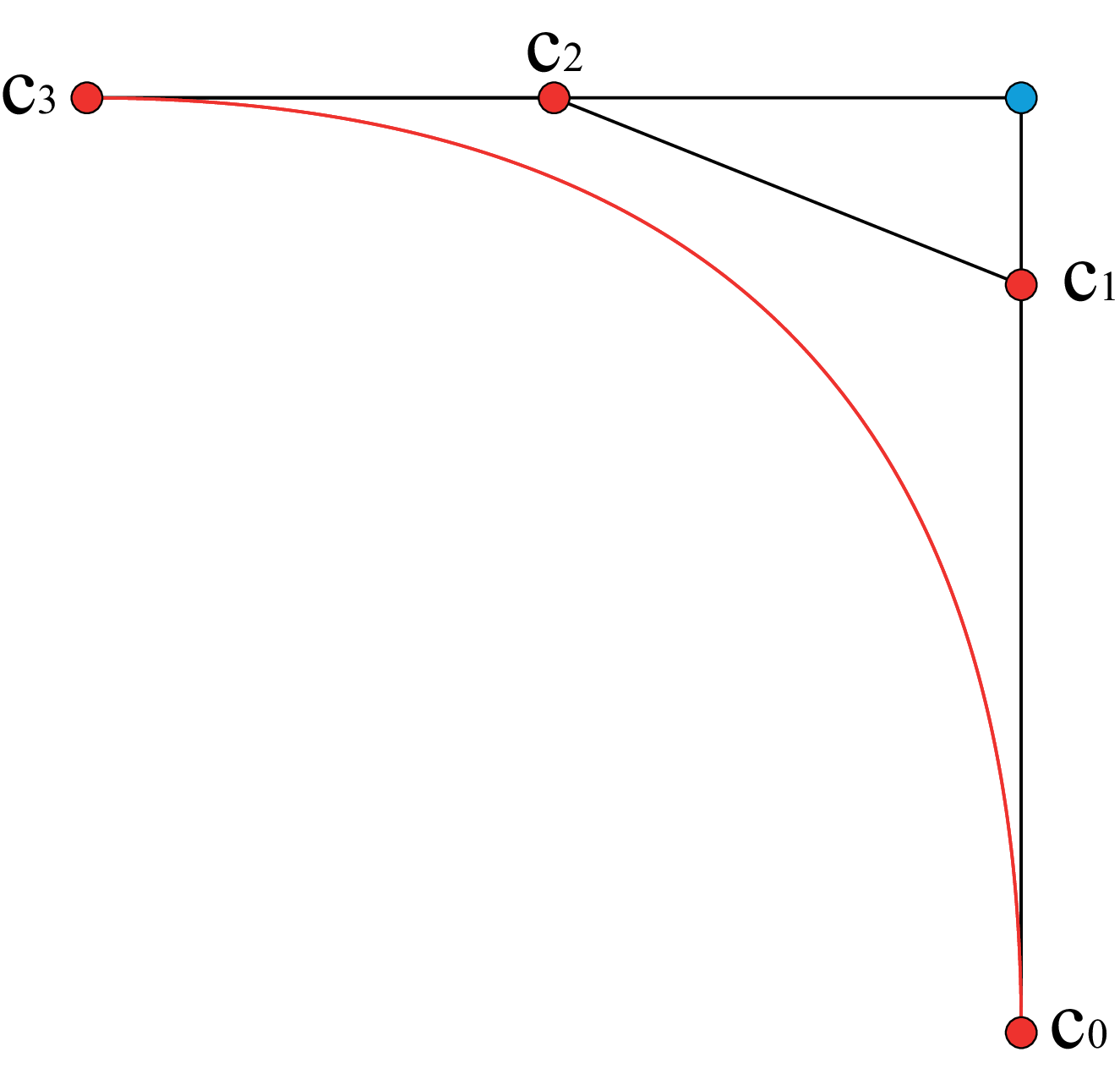}
\caption{Arc of parabola with control polygon $\{(1,0),(1,1),(0,1)\}$
seen as a rational cubic with control polygon  $\{(1,0), (1, 4/5),
(1/2, 1), (0, 1)\}$ and weights $\{2, 5/3, 4/3,1\}$}
\label{fake1}
\end{figure}
We have already seen in Example~\ref{fakecubic} that the cubic parametrisation is reducible and 
describes in fact an arc of conic.

In order to obtain an irreducible parametrisation, we calculate the 
gcd of the numerator and the denominator of the parametrisation.

The complex control polygon for the curve is $\{1, 1+4i/5,1/2+i,i\}$ 
and for the numerator $p(t)$ we get $\{2, 5/3+4i/3,2/3+4i/3,i\}$, $\lbrac2, 
5+4i,2+4i,i\rbrac$ as reduced coefficients. For the denominator $q(t)$ 
we have $\lbrac2, 5, 4,1\rbrac$ as reduced coefficients.

The quotient of $p(t)$ and $q(t)$ is $i$ and the remainder $r_{1}(t)$
has reduced coefficients $\lbrac2 - 2i, 5 - i, 2\rbrac$.

As next step of Euclid's algorithm, we divide $q(t)$ by $r_{1}(t)$ 
and we get a remainder $r_{2}(t)$ with coefficients $\{i,i/2\}$.

We divide $r_{1}(t)$ by $r_{2}(t)$ 
and get a null remainder and hence 
$\mathrm{gcd}(p(t),q(t))=r_{2}(t)$.

Dividing $q(t)$ by $r_{2}(t)$,
we get as list of reduced weights $\lbrac-2i,-4i,-2i\rbrac$ and 
$\{-2i,-2i,-2i\}$ as list of complex weights. Since all weights are equal, it is in fact an arc 
of parabola.

For the quotient of $p(t)$ and $r_{2}(t)$,
we get  for the numerator $\lbrac-2i, 4 - 4i, 2\rbrac$ as reduced  and 
$\{-2i, 2 - 2i, 2\}$ as coefficients. After factoring the weights we get $\{1,1+i, 
i\}$ as complex control polygon for the arc of parabola.

Hence our cubic arc is an arc of parabola with control polygon 
$\{(1,0),(1,1),(0,1)\}$. We see both control polygons in 
Figure~\ref{fake1}. \end{example}

\section{From complex to real\label{real}}

We have seen how to get complex parametrisations starting with real
ones.  To go from complex parametrisations to real ones, it is simple.
If we have a complex rational parametrisation of the form
$c(t)=p(t)/q(t)$, we may produce a real parametrisation in the usual
way, multiplying by the complex conjugate of the denominator
$\overline{q(t)}$, $c(t)=p(t)\overline{q(t)}/q(t)\overline{q(t)}$,
which doubles formally the degree of the parametrisation.

If the complex parametrisation $c(t)=p(t)/q(t)$ is described by a
control polygon $\{z_{0},\ldots,z_{n}\}$ and a list of complex weights
$\{\omega_{0},\ldots,\omega_{n}\}$, then the coefficients for the
numerator $p(t)$ are $\{\omega_{0}z_{0},\ldots,\omega_{n}z_{n}\}$.  The
control polygon for $\overline{q(t)}$ is simply
$\{\overline{\omega_{0}},\ldots,\overline{\omega_{n}}\}$.

According to (\ref{multiply}), multiplication by $\overline{q(t)}$ is 
accomplished using reduced control coefficients. For the denominator we 
have a new list of reduced weights,
\[\hspace{-0.5cm}\{\tilde w_{0},\ldots,\tilde w_{2n}\}, \quad
\tilde w_{I}=\sum_{j+k=I}\tilde \omega_{j} 
\overline{\tilde\omega_{k}},\ I=0,\ldots,2n,\]
and a new reduced control polygon for the numerator,
\[\hspace{-0.5cm}\{\tilde p_{0},\ldots,\tilde p_{2n}\}, \quad
\tilde p_{I}=\sum_{j+k=I}z_{j}\tilde \omega_{j} 
\overline{\tilde\omega_{k}},\ I=0,\ldots,2n.\]

Let's see it with an example.  \begin{example}Arc of a circle defined
by a control polygon $\{1,i\}$ and a list of weights,
$\{1+i,2\}$.

We have already seen this parametrisation in Example~\ref{excircle}. 
For a real parametrisation of degree two, we have a list of weights,
\[\lbrac1+i,2\rbrac*\lbrac1-i,2\rbrac=\lbrac2,4,4\rbrac,\]
and for the numerator, with control polygon $\{1+i,2i\}$ we get
\[\lbrac1+i,2i\rbrac*\lbrac1-i,2\rbrac=\lbrac2,4+4i,4i\rbrac,\]
and after factoring the weights we get as control 
polygon, $\{1,1+i,i\}$.

That is, the real control polygon for the arc is 
$\{(1,0),(1,1),(0,1)\}$ and the list of weights is $\{2,2,4\}$ (or 
$\{1,1,2\})$) as we already knew.\end{example}

\section{Some examples of curves\label{curves}}

There are a number of curves in geometry which can be constructed
using circles, such as trochoids \cite{reyes-complex}, which are
amenable for complex representation.  Using inversion of plane curves,
there is a number of other classical constructions which take
advantage of this framework \cite{walker}.  A catalog of cases may be
found in \cite{lawrence}.

Geometric inversion with respect to a circle of radius $R$ centered at
the origin is defined as $R^{2}/\bar{z}$.  Since most of the examples
are symmetric, we shall use the complex inversion $1/z$ instead, which 
amounts of considering $R=1$ and performing a reflection with respect 
to the real axis. The effect of complex inversion on complex control 
polygons and weights has already been discussed in 
Section~\ref{define}.

It can be verified that the geometric inversion of a conic with
respect to one of its points gives a circular rational cubic, for
example, the cissoid of Diocles, while when the center of inversion
does not belong to the conic, the inverted curve is a bicircular
rational quartic, for instance, the cardioid and the Bernouilli
lemniscate  \cite{lawrence}.

%
%

In the following examples we draw the original curve in blue, the 
inverted curve in red, showing with a dashed line the parts which are 
not included in the parametrisation.

\subsection{Cissoid of Diocles}

This curve has been used for solving the problem of doubling a cube:
given a cube, find another cube whose volume is double that of the
first one.  

A rational parametrisation of the cissoid, expressed in complex form, for constant
$a\in\mathbb{R}\setminus{\{0\}}$, is
\[\frac{2a s^{2}(s-i)}{s^{2}+1}=\frac{2a s^{2}}{s+i},\quad 
s\in\mathbb{R}.\]
Some terms can be factored and the degree becomes 
two.

This parametrisation is the inverse of the parametrisation of a 
parabola $\frac{s+i}{2a s^{2}}$ with equation $y=2ax^{2}$. 
Another way of constructing a cissoid is inverting a parabola 
with respect to its vertex.
\begin{figure}
\centering
\includegraphics[height=1.5cm]{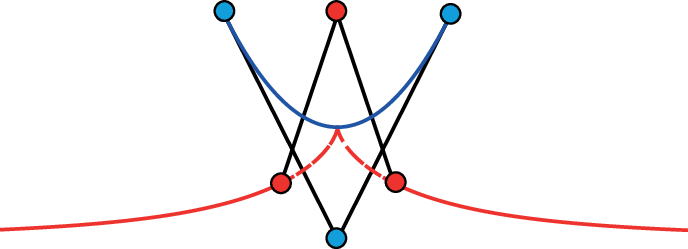}
\caption{Arc of a cissoid of Diocles}
\label{cissoid}
\end{figure}

Taking $a=1/2$, if we start with an arc of the parabola of equation
$y=x^{2}$, with vertex at the origin, with complex control polygon
$\{-1+i,-i,1+i\}$ and trivial weights, after complex inversion we
obtain an arc of cissoid with complex control polygon $\{-1/2 - i/2,
i, 1/2 - i/2\}$ and complex weights, $\{-1 + i, -i, 1 + i\}$.

Since the arc of parabola comprises its vertex, due to inversion, the arc of cissoid 
extends to infinity (see Figure~\ref{cissoid}).

\subsection{Cardioid}

A cardioid is an example of trochoid with a single cusp 
\cite{reyes-complex}. A complex rational parametrisation can be
\[2a\frac{1-s^2 +2is}{(s^2+1)^2}=-\frac{2a}{(s+i)^{2}},\quad 
s\in\mathbb{R},\]
for constant $a\in\mathbb{R}\setminus{\{0\}}$, which allows us to lower the degree from four to two.

The inverse of this curve is parametrised as $-(s+i)^{2}/2a$ and it is hence a parabola with equation 
$x=(1-a^{2}y^{2})/2a$.

\begin{figure}
\centering
\includegraphics[height=4cm]{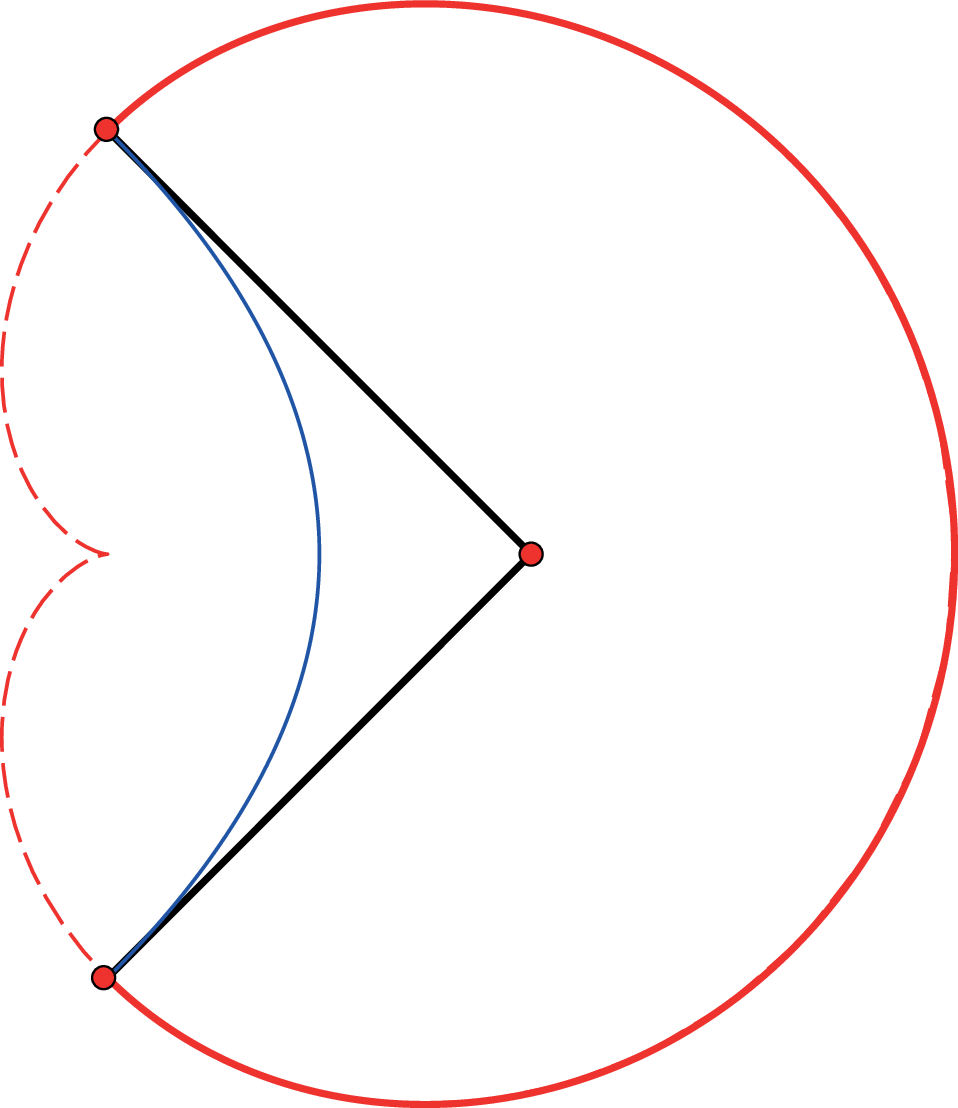}
\caption{Arc of a cardioid as inverse of a parabola}
\label{cardioid}
\end{figure}

A cardioid can be constructed as inverse of a parabola with respect to
its focus.  For instance, take $a=1$ and consider the parabola with
equation $x=(1-y^{2})/2$ and focus at the origin.  An arc of this
parabola is described by the complex control polygon $\{-i, 1, i\}$ and
trivial weights.  By complex inversion we get an arc of cardioid with
complex control polygon $\{i, 1,- i\}$ and a list of complex weights
$\{-i, 1, i\}$ (see Figure~\ref{cardioid}).

\subsection{Lemniscate of Bernoulli}

Another well known planar curve is the lemniscate of Bernoulli.  This
curve may be seen as the inverse of an equilateral hyperbola with
respect to its center.  For a hyperbola of equation
$x^{2}-y^{2}=a^{2}$, centered at the origin, a complex rational
parametrisation for the lemniscate can be
\[a\frac{s+s^{3}+ i (s-s^{3})}{1+s^{4}}=\frac{a(1 -i)s}{s^2 -i}, \quad 
s\in\mathbb{R},\]
for constant $a\in\mathbb{R}\setminus{\{0\}}$.

For simplicity, we take $a=1$ and an arc of hyperbola defined by a 
complex control polygon, $\{\sqrt{2} - i, \sqrt{2}/{2}, 
\sqrt{2} + i\}$ with weights $\{1,\sqrt{2},1\}$.

\begin{figure}
\centering
\includegraphics[height=4cm]{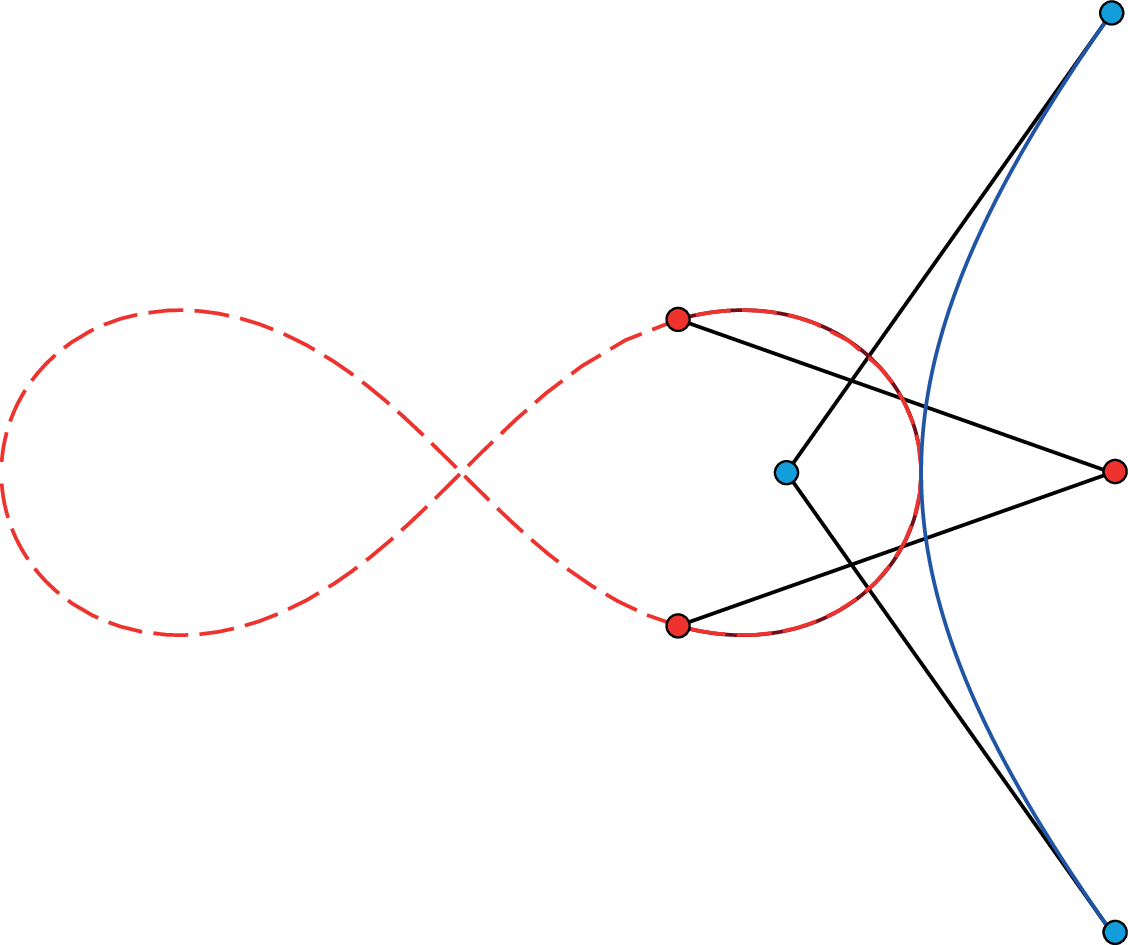}
\caption{Arc of a lemniscate (red) as inverse of a arc of an equilateral 
hyperbola (blue)}
\label{lemniscate}
\end{figure}

The complex control polygon for the arc of lemniscate is 
$\{\frac{\sqrt{2}+i}{3},\sqrt{2},\frac{\sqrt{2}-i}{3}\}$ and the list of weights, 
$\{\sqrt{2}-i, 1, \sqrt{2}+i\}$. 
In Figure~\ref{lemniscate} we show both arcs in red and blue.

%

%
%
%
%
%
%

\section{Conclusions\label{conclude}}

In this paper we have developed the complex framework for rational
B\'ezier plane curves, showing its advantages in terms of degree
lowering and the possibility of using two different groups of
projective transformations for manipulating plane rational curves,
while preserving their properties, going further than with the
standard real framework.  In fact, the complex framework can be
considered as already included in the CAD paradigm, since it uses the
same elements, vertices and weigths, and constructions and just
involves extending the type of variable for weights from real to
complex.  

We have shown how to change from one framework to another, resorting 
to degree-elevation and factorisation of polynomials, making using of 
classical algebraic elements, such as resultants and Euclid's 
theorem, but for polynonials in Bernstein basis.

In particular, we have shown a simple procedure to check if a
parametrisation is reducible or not in terms of the control polygon
and weights.

Examples of curves for which the complex framework enables the 
possibility of using lower degree parametrisations are 
included. 

\section*{Acknowledgments}

This work is partially supported by the Spanish Agencia Estatal de 
Investigaci\'on through research grant PID2024-158664NB-C21.

\bibliographystyle{elsarticle-num-names}
\bibliography{cagd}

\end{document}